\documentclass[11pt]{amsart}
\usepackage{amssymb}

\usepackage{palatino}
\input amssym.def

\usepackage{amsmath, amsfonts}
\usepackage{amssymb}
\usepackage{amscd}
\usepackage[mathscr]{eucal}
\usepackage{palatino}
\newfont{\cyrr}{wncyr10}

\newcommand{\Q}{{\mathbb Q}}

\newcommand{\N}{{\mathbb N}}
\newcommand{\bQ}{\overline{\Q}}

\renewcommand{\mod}{{\, \rm mod \, }}

\newtheorem{thm}{Theorem}

\newtheorem{rmk}{Remark}[section] 

\newtheorem{conj}[thm]{Conjecture}
\newcommand{\thmref}[1]{Theorem~\ref{#1}}

\begin{document}

\title{Linear and algebraic independence of Generalized 
Euler-Briggs constants}

\author{Sanoli Gun, V. Kumar Murty and Ekata Saha}
\address{Sanoli Gun and Ekata Saha\\ \newline
Institute of Mathematical Sciences, C.I.T. Campus, Taramani, 
Chennai, 600 113, India.}
\email{sanoli@imsc.res.in}  
\email{ekatas@imsc.res.in}

\address{V. Kumar Murty\\ \newline
Department of Mathematics, University of Toronto,
40 St. George Street, Toronto, ON, Canada, M5S 2E4.} 
\email{murty@math.toronto.edu}

\subjclass[2010]{11J81, 11J91}

\keywords{Generalized Euler-Briggs constants, Baker's theory of 
linear forms in logarithms, Weak Schanuel's conjecture}

\maketitle

\begin{abstract}
Possible transcendental nature of Euler's constant $\gamma$
has been the focus of study for sometime now. 
One possible approach  is to consider $\gamma$ not in isolation, but 
as an element  of the  infinite family of generalised 
Euler-Briggs constants. In a recent work \cite{GSS}, 
it is  shown that the infinite list of  generalized
Euler-Briggs constants  can have at most one algebraic number. 
In this paper, we study the dimension
of spaces generated by these generalized 
Euler-Briggs constants over number fields.
More precisely, we obtain non-trivial lower
bounds (see \thmref{pre} and \thmref{linear-ind}) on the 
dimension of these
spaces and consequently establish the infinite dimensionality
of the space spanned.
Further, we study linear and algebraic
independence of these constants over the field
of all algebraic numbers.
\end{abstract}

\section{\large Introduction}

In 1731, Euler introduced the following constant 
$$
\gamma:=
\lim_{ x \to \infty} (\sum_{n \le x}
\frac{1}{n}  ~-~ \log x )
$$ 
and derived a number of identities involving  
$\gamma$, special values of the Riemann zeta function
and other known constants. 
After Euler several other notable mathematicians 
including Gauss and Ramanujan have studied this 
constant in depth.
For a beautiful account of the various aspects
of research about this constant, we refer the
reader to a recent article of Lagarias~\cite{JL}.

The appearance of $\gamma$
in its various avatars makes it a fundamental 
object of study in number theory. 
Though we expect that $\gamma$ is transcendental, 
it is not even known to be irrational.
However there are some transcendence results
involving $\gamma$. To the best of our knowledge,
the first such result was due to Mahler \cite{KM}.
He proved that for any non-zero algebraic
number $\alpha$, the number
$$
\frac{\pi Y_0(\alpha)}{2J_0(\alpha)} 
- \log \frac{\alpha}{2} - \gamma
$$
is transcendental, where $J_0$ and $Y_0$
are Bessel functions of the first and second kind 
of order zero. More precisely, 
\begin{eqnarray*}
J_0(\alpha) 
& :=& 
\sum_{n=0}^{\infty} \frac{(-1)^n}{(n!)^2}(\frac{\alpha}{2})^{2n}~, 
\phantom{m}
H_n := \sum_{j=1}^n \frac{1}{j} \\
\text{and}\phantom{m}  
\frac{\pi}{2} Y_0(\alpha) 
&:=& 
\left( \log(\frac{\alpha}{2}) + \gamma \right) 
J_0(\alpha)  ~+~ 
\sum_{n =1}^{\infty} (-1)^{n-1}
 \frac{H_n}{(n!)^2}(\frac{\alpha^2}{4})^n. 
\end{eqnarray*}

In the rather difficult subject of transcendence,
sometimes it is  more pragmatic to look  at a  
family of numbers as opposed to a single number and 
derive something meaningful. There are two significant results 
which are worth mentioning at this point.
The first one is by R. Murty and Saradha \cite{MS2}. 
They proved the following theorem.
\begin{thm}{\rm (Murty and Saradha)}~
Let $q> 1$ be a natural number. 
At most one of the numbers 
$$
\gamma, ~ \gamma(a, q),  1 \le a \le q, ~(a, q) =1
$$
is algebraic. Here 
$$
\gamma(a, q):= \lim_{x \to \infty} 
( \sum_{n \le x \atop n \equiv a \mod{q}} \frac{1}{n}  
-  \frac{1}{q}\log x).
$$
\end{thm}
The constants $\gamma(a,q)$ were introduced by Briggs \cite{WB}
and later studied extensively by Lehmer \cite{DL}. We will call
these constants as Euler-Briggs constants. 
The second result involving the transcendence of 
$\gamma$ is due to Rivoal \cite{TR} and  Pilehrood-Pilehrood \cite {PP} who
proved the following theorem.
\begin{thm}{\rm (Rivoal / Pilehrood-Pilehrood)}~
At least one of the numbers $\gamma$ and 
$\delta := \int_{0}^{\infty} \frac{e^{-w}}{1+w} dw$
is transcendental. 
\end{thm}
The constant $\delta$ is known as the Euler-Gompertz constant.
Note that the constants $\gamma$ and $\delta$
are part of a family of numbers called 
exponential periods (see \cite{MK}
and also page 595 of \cite{JL}) introduced by
Kontsevich and Zagier \cite{KZ}.

Another important set of numbers with
$\gamma$ as a member were introduced
by Diamond and Ford. In 2008, while
studying the  Riemann hypothesis, they
introduced the so-called generalized Euler's constants. 
In order to introduce these numbers, let
us set some notation.

Throughout the paper,
$P$ will denote the set of prime numbers
and $p$ will denote a prime number.
For any finite subset $\Omega \subset P$, set
\begin{eqnarray}\label{imp}
P_{\Omega} 
&:=& \begin{cases}
\prod_{p \in \Omega} p & \mbox{ if   }~ \Omega \ne \phi, \\
1  & \mbox{ otherwise},
\end{cases}  \\
\phantom{m}
\text{ and }
\phantom{mm}
\delta_\Omega 
&:=& 
\begin{cases}
\prod_{p \in \Omega} ( 1- \frac{1}{p})
& \mbox{ if   }~ \Omega \ne \phi,  \nonumber \\
1  & \mbox{ otherwise}.
\end{cases} 
\end{eqnarray}
For any finite subset $\Omega \subset P$, 
Diamond and Ford \cite{DF} defined the
generalized Euler's constant as follows;
$$
\gamma(\Omega) 
: = 
\lim_{ x \to \infty} (\sum_{n \le x \atop (n, P_{\Omega})=1}
\frac{1}{n}  ~-~ \delta_{\Omega} \log x ).
$$
Note that, when $\Omega = \phi$, 
then $\gamma(\Omega) = \gamma$.
In a recent work, R. Murty and Zaytseva \cite{MZ}
noticed that at most  one number in the infinite list
$\gamma(\Omega)$ as $\Omega$ varies over
finite subset of primes is algebraic.

Following Euler-Briggs and Diamond-Ford, one can 
now define for any finite subset $\Omega \subset P$ and 
natural numbers $a, q$ 
with $(q, P_{\Omega})=1$, the constants
$$
\gamma(\Omega, a, q)
~:=~ 
\lim_{ x \to \infty} 
(\sum_{\substack{ n \le x \atop {n \equiv a \mod q \atop (n, P_{\Omega})=1}}}
\frac{1}{n} 
   ~~-~   \delta_{\Omega} \frac{\log x}{q} ~).
$$  
When $\Omega = \phi$, we see that
$\gamma(\Omega, a , q) = \gamma(a,q)$, the classical 
Euler- Briggs constants. 
From now on, we refer to the constants
$\gamma(\Omega, a, q)$ as 
generalized Euler-Briggs constants.    
Moreover, when $q=1$ and
$\Omega \subset P$ is a finite set, we have
$$
\gamma(\Omega, a, 1) 
= 
\gamma(\Omega) 
~=~ 
\delta_{\Omega} ~(\gamma ~+~ \sum_{p \in \Omega}\frac{\log p}{p-1}),
\text{   where  } a \in \N.
$$
The last equality has been established by Diamond and 
Ford in \cite{DF}.
In a recent work, the first and the third author along with 
Sneh Bala Sinha \cite{GSS} (see also \cite{GSS1})
proved the following results;
\begin{thm}{\rm (Gun, Saha and Sinha)}\label{one}
Let $a$ and $q> 1$ be natural numbers with $(a,q)=1$ and $S$
be the set of prime divisors of~$q$.
Also let 
$$
U:= \left\{~\Omega \mid \Omega \text{ is a finite set of primes},  
~\Omega \cap S = \phi\right\}.
$$
Then the set $T := \left\{ \gamma(\Omega, a, q) \mid \Omega \in U \right\}$
is infinite and has at most one algebraic element. 
\end{thm}

\begin{thm}\label{two}{\rm (Gun, Saha and Sinha)}~
Let $\Omega$ be a finite set of primes and 
$S= \{ q_1, q_2, \cdots \}$ be an infinite set of mutually co-prime
natural numbers $q_i >1$ with $(q_i, P_{\Omega}) =1$
for all $i \in \N$.  Then for any $a \in \N$ 
with $(a, q_i)=1 ~\forall i$, the set
$$
T:=  \{ ~\gamma(\Omega, a, q_i)  ~\mid~  q_i \in S \}
$$
has  at most one algebraic element.   
\end{thm}

In order to prove these theorems, one needs to find
a closed formula for generalized Euler-Briggs constants.
In \cite{GSS}, it was proved that
\begin{eqnarray}\label{br}
\gamma(\Omega, a, q)  
&=& \frac{1}{\varphi(q)}
\sum_{\substack{\chi \mod q \atop \chi \ne \chi_0}}
\overline{\chi}(a) L(1,\chi)
\prod_{p\in\Omega} (1-\frac{\chi(p)}{p})
~+~ \frac{\delta_{\Omega}}{q} (\gamma + \sum_{p \mid q} 
\frac{\log p}{p-1} +
\sum_{p \in \Omega}\frac{\log p}{p-1}),
\end{eqnarray}
where $\Omega$ is a finite subset of primes, $a,q$ are co-prime
natural numbers with $(q, P_{\Omega})=1$.

The articles \cite{{GSS, GSS1}} can be thought of a 
generalization of the work of R. Murty and Zaytseva \cite{MZ} on 
generalized Euler constants. Whereas the 
results of R. Murty and Zaytseva are obtained by 
Hermite and Lindemann theorem, 
the proofs in \cite{GSS} require careful analysis of 
units in cyclotomic fields
and Baker's theorem on linear forms in logarithms. 

The above theorems do not answer the question
of linear independence of these constants
over a number field or over $\bQ$.

In this article, we establish non-trivial  
lower bounds for the dimension of
these spaces. To start with,
we have the following theorem over
the field of rational numbers $\Q$.

\begin{thm}\label{pre}
Let $\Omega \subset P$ be a finite subset of primes
and $P_{\Omega}$ be as in \eqref{imp}. 
Consider the $\Q$-vector space
$$
V_{\Q,N} := \Q \left <  \gamma(\Omega, m, n) ~|~   
1\le m\le n \le N,~ ~(m, n)=1=(n, P_{\Omega}) \right >.
$$
Then for $N$ sufficiently large, we have
$$
N ~\ll_\Omega~~ \dim_{\Q} V_{\Q,N},
$$
where the implied constant depend
on $\Omega$. In particular, the dimension of the
$\Q$-vector space
$$
V_{\Q} := \Q \left <  \gamma(\Omega, m, n) ~|~   
m, n \in \N, ~(m, n)=1=(n, P_{\Omega}) \right >
$$
is infinite.
\end{thm}

In fact, one has the following general theorem
about linear independence of these constants over
number fields.  
\begin{thm}\label{linear-ind}
Let $K$ be a number field with discriminant $d > 1$, 
$\Omega \subset P$ be a finite
subset of primes, $P_{\Omega}$ be as in \eqref{imp}
such that $K \cap \Q(\zeta_{P_{\Omega}}) = \Q$, where 
$\zeta_{P_{\Omega}} := e^{\frac{2\pi i}{P_{\Omega}}}$.
Consider the $K$-vector space
$$
V_{K,N} := K \left <  \gamma(\Omega, m, n) ~|~   
1\le m\le n \le N, ~~(m, n)=1=(n, dP_{\Omega}) \right >.
$$
Then for $N$ sufficiently large, we have
$$
N ~\ll_{K, \Omega}~~  \dim_{K} V_{K,N},
$$
where the implied constant depend on $\Omega$ and $K$.
In particular, the $K$-vector space
$$
V_{K} := K \left <  \gamma(\Omega, m, n) ~|~   
m, n \in \N, ~(m, n)=1=(n, dP_{\Omega}) \right >
$$
is infinite dimensional.
\end{thm}

\begin{rmk}
Note the trivial upper bounds for
$\dim_{\Q} V_{\Q,N}$ in \thmref{pre} and for $\dim_{K} V_{K,N}$ in
\thmref{linear-ind} are $N^2$.
\end{rmk}

Next we study the linear independence of these constants
over the field of algebraic numbers.  In order to do so,
let us set
\begin{eqnarray*}
C(q) &:=& \{ \Omega \subset P ~|~~  |\Omega| < \infty, 
~(q, P_{\Omega})=1 \}, \phantom{m} \text{ where } q\in \N.
\end{eqnarray*}
We define an equivalence relation on the set of all 
$\gamma(\Omega, a, q)$'s as $\Omega$ varies over elements of $C(q)$
and $a,q$ are co-prime natural numbers. 
We say that $\gamma(\Omega_1, a, q)$ and $\gamma(\Omega_2, a, q)$
are equivalent, denoted by $\gamma(\Omega_1, a, q) 
\sim \gamma(\Omega_2, a, q)$, 
if there exists $\lambda \in {\bQ} \setminus \{ 0\}$ such that 
$\gamma(\Omega_1, a, q) = \lambda \gamma(\Omega_2, a, q)$. 

In this set-up, we prove the following theorem.
\begin{thm}\label{E2}
Let $a, q$ be natural numbers with $(a, q)=1$.
Consider the set
$$
M_1 := \{ \gamma(\Omega, a , q) ~|~ \Omega \in C(q) \}.
$$
Then each equivalence class $[\gamma(\Omega, a, q)]$ in $M_1$ 
has at most two elements. 
\end{thm}

Next let $\Omega \subset P$ be a finite set, $P_\Omega$
be as in \eqref{imp},  $a \in \N$ 
and 
$$
C(a, \Omega) := \{ q \in \N  
~|~~  (a,q)=1 = (q, P_{\Omega}) \}.
$$
As before, one can define an equivalence relation 
on the set $\gamma(\Omega, a, q)$'s, where
$q \in C (a, \Omega)$.
In this set-up, we prove;

\begin{thm}\label{E3}
Let $\Omega$ be a finite set of primes, $\{q_i\}$ be a sequence of 
mutually co-prime natural numbers and $a$ be a natural number
such that $(a, q_i)=1$ for all $i$. Consider the set
$$
M_2 := \{ \gamma(\Omega, a , q_i) ~|~ q_i \in C (a, \Omega) \}.
$$
Then each equivalence class $[\gamma(\Omega, a, q_i)]$ in $M_2$ 
has at most two elements. 
\end{thm}
We see that \thmref{E2} and \thmref{E3}
give information about pairwise $\bQ$-independence 
of the generalized Euler-Briggs constants.
However, they do not say that the vector
space generated by these constants 
over $\bQ$ is infinite dimensional. 
In this regard, we have the 
following theorem;

\begin{thm}\label{closure}
Let $a,q$ be natural numbers with $(a,q) =1$.
Then the dimension of the $\bQ$-vector space
$$
V_{\bQ} : = \bQ < \gamma(\Omega, a, q)  ~|~ \Omega \in C(q)>
$$
is infinite over $\bQ$.
\end{thm}

We end this section with a brief outline of the
structure of the paper. In \S2, we list the
inputs from transcendence theory
relevant to our work. We also state a general
non-vanishing result (\thmref{ESS2})
which is integral to our investigation. The proof
of this theorem is detailed in \S3. 
We devote \S4 to prove all the linear independence
results indicated in the introduction. Finally, in
\S5, we state and prove some algebraic independence
results assuming the Weak Schanuel conjecture.

\smallskip 
 
\section{Preliminaries}

\smallskip 

In this section, we state the theorems which will be required
to prove our results. The first and the third author proved 
the following theorem about the existence
of an infinite sum (see~\cite{GS}).

\begin{thm}\label{gs}{\rm (Gun and Saha)}~
Let $f$ be a periodic arithmetic function
with period $q \ge 1$ and $M$ be a natural number
co-prime to $q$. Then 
$$
\sum_{n \ge 1 \atop (n, M)=1} \frac{f(n)}{n}
$$ 
exists if and only if $~\sum_{a=1}^q f(a) =0$.
Moreover, whenever the above sum exists, we have
$$
\sum_{n \ge 1 \atop (n, M)=1} \frac{f(n)}{n}
=
\sum_{a=1}^q f(a) \gamma(\Omega, a, q),
$$
where $\Omega$ is the set of prime divisors of $M$.
\end{thm}
An important ingredient to prove
\thmref{linear-ind}
is the following non-vanishing 
result of Baker, Birch and Wirsing
(see \cite{BBW}, see also chapter 23 of \cite{MR}).
\begin{thm}{\rm (Baker, Birch and Wirsing).}\label{non-vanish}
Let $f$ be a non-zero algebraic valued 
periodic function with period $q$
defined on the set of integers.
Also let $f(n)=0$ whenever $1 < (n,q) < q$ and
the $q$-th cyclotomic polynomial $\Phi_q(X)$ 
be irreducible over $\Q(f(1),\cdots,f(q))$, then
$$
\displaystyle \sum_{n=1}^{\infty} \frac{f(n)}{n} ~~\neq 0.
$$
\end{thm}
Other important theorems that are required to prove 
our results are the following.
\begin{thm}~\label{ESS2}
Let $q_1, q_2 , q_3 > 1$ be mutually co-prime natural numbers. 
Then for any algebraic numbers $\alpha_p, \beta_{\chi}, \beta_{\phi}, \beta_{\psi}$,
the number
$$
\sum_{p | q_1q_2q_3} \alpha_p \log p  
~+~ \sum_{\chi \mod {q_1} \atop \chi \ne \chi_0} \beta_{\chi} L(1, \chi)
~+~ \sum_{\phi \mod {q_2} \atop \phi \ne \phi_0} \beta_{\phi} L(1, \phi)
~+~ \sum_{\psi \mod {q_3} \atop \psi \ne \psi_0} \beta_{\chi} L(1, \psi)
$$
is transcendental provided not all $\alpha_p, \beta_{\chi}, \beta_{\phi}, \beta_{\psi}$
for even characters $\chi, \phi, \psi$ are zero.
\end{thm}
This result is new and we will give a proof of this result
in the next section. A particular
case of \thmref{ESS2} was noticed in \cite{GSS}.
In the same paper, the authors also proved the following theorem
which will be required to prove \thmref{closure}.
\begin{thm}\label{cl}{\rm (Gun, Saha and Sinha)}~
Let $q > 1$ be a natural number and $\Omega_1, 
\cdots, \Omega_t \in C(q)$ be disjoint subsets of
prime numbers. Then for any algebraic numbers $\alpha_p, \beta_{\chi},
\epsilon_{\Omega_i, p}$, the number
$$
\sum_{p|q} \alpha_p \log p  
~+~ \sum_{i=1}^t \sum_{p \in \Omega_i} \epsilon_{\Omega_i, p} \log p
~+~ \sum_{\chi \mod {q} \atop \chi \ne \chi_0} \beta_{\chi} L(1, \chi)
$$
is transcendental provided not all $\alpha_p, \epsilon_{\Omega_i, p}, 
\beta_{\chi}$ for even characters $\chi$ are zero.
\end{thm}

We shall be using the following 
result of Baker (see pages 10 and 11 of \cite{AB},
see also chapter 19 of \cite{MR}). 
\begin{thm} {\rm(Baker)} \label{B}
Let $\alpha_1, \cdots, \alpha_n \in \bQ \setminus \{0\}$ and 
$\beta_1, \cdots, \beta_n \in \bQ$, then
$$
\beta_1 \log \alpha_1  + \cdots + \beta_n \log \alpha_n
$$ 
is either zero or transcendental. The latter case 
arises if $\log \alpha_1, \cdots, \log \alpha_n$
are linearly independent over~$\Q$ and not all
$\beta_1, \cdots, \beta_n$ are zero.
\end{thm}
Finally, we will be using the following result 
about the non-vanishing of certain special linear forms in logarithms of 
non-zero algebraic numbers.
\begin{thm}{\rm(R. Murty and Saradha \cite{MS2}, see 
also R. Murty and K. Murty \cite{MM})}\label{pione}
Let $\alpha_1, \cdots, \alpha_n$ be positive algebraic 
numbers. If $\beta_0, \cdots, \beta_n$
are algebraic numbers with $\beta_0 \ne 0$, then
$$
\beta_0 \pi  ~~+~ \sum_{i=1}^{n} \beta_i \log\alpha_i 
$$
is a transcendental number and hence non-zero.
\end{thm}

\section{Proof of \thmref{ESS2}}

We will prove this theorem by contradiction.
We know that for any even Dirichlet character 
$\chi \ne \chi_0$, one can write
$L(1, \chi)$ as a non-zero algebraic multiple of 
\begin{equation}\label{key}
\sum_{1 < a < q/2 \atop (a,q)=1} \overline{\chi}(a) \log \xi_a,
\end{equation}
where $\xi_a$'s are real multiplicatively independent units in the
cyclotomic field $\Q(\zeta_q)$, known as Ramachandra units  
(see pages 147 to 149 of \cite{KR}, page 149 of \cite{LW} 
as well as page 1728 of \cite{MM}).  
For any odd Dirichlet character~$\chi$, we know that
$L(1,\chi)$ is a non-zero algebraic multiple of $\pi$ 
(see page 38 of \cite{LW}). 
Using these results and \thmref{pione}, 
we can therefore ignore the odd characters. In order to complete
the proof of the theorem, we will now show that 
\begin{enumerate}
\item
$\log p$~:~~ for all primes $p| q_1q_2q_3$  
\item
$L(1,\chi)$~:~~ for all  even non-principal characters $\chi$ modulo $q_1$ 
\item
$L(1,\phi)$~:~~ for all even non-principal characters $\phi$ modulo $q_2$
\item
$L(1,\psi)$~:~~ for all even non-principal characters $\psi$ modulo $q_3$
\end{enumerate}
are linearly independent over $\bQ$. 
Suppose not.  Then there exists algebraic numbers $\alpha_p$
for $p|q_1q_2q_3$ and $\beta_{\chi}, \beta_{\phi}, \beta_{\psi}$, where 
$\chi, \phi, \psi$ vary over non-principal even Dirichlet characters
modulo $q_1$, $q_2$ and $q_3$ respectively, not all zero, such that
$$
\sum_{p | q_1q_2q_3} \alpha_p \log p 
~+~  \sum_{\chi \text{ even } \atop { \chi \ne \chi_0}} \beta_{\chi} L(1, \chi)
~+~  \sum_{\phi \text{ even } \atop { \phi \ne \phi_0}} \beta_{\phi} L(1, \phi)
~+~  \sum_{\psi \text{ even } \atop { \psi \ne \psi_0}} \beta_{\psi} L(1, \psi)
~=~  0.
$$
We can rewrite the above expression as
$$
\sum_{p | q_1q_2q_3} \alpha_p \log p
~+~  \sum_{1< a< q_1/2 \atop (a, q_1) =1}  \delta_a \log \xi_a 
~+~   \sum_{1< b< q_2/2 \atop (b, q_2) =1}  \delta_b \log \xi_b  
~+~   \sum_{1< c< q_3/2 \atop (c, q_3) =1}  \delta_c \log \xi_c 
~=~ 0, 
$$
where $\xi_a, \xi_b, \xi_c$ 's are multiplicatively independent units
in $\Q(\zeta_{q_1}), \Q(\zeta_{q_2}) $ and $\Q(\zeta_{q_3})$ respectively.
Now by Baker's Theorem, we have
\begin{equation}\label{extra1}
\prod_{p | q_1q_2} p^{c_p}  
~=~ 
\prod_{1 < a < q_1/2 \atop (a, q_1) =1} {\xi}_a^{d_a} 
\prod_{1 < b < q_2/2 \atop (b, q_2) =1} {\xi}_b^{e_b}
\prod_{1 < c < q_3/2 \atop (c, q_3) =1} {\xi}_c^{f_c}
\end{equation}
where $c_p, d_a, e_b, f_c$'s are integers.
By taking norms on both sides of \eqref{extra1}, we get $c_p =0$ 
for all $p$. Hence
\begin{equation}\label{re}
\prod_{1 < a < q_1/2 \atop (a, q_1) =1} {\xi}_a^{d_a} 
~=~ 
\prod_{1 < b < q_2/2 \atop (b, q_2) =1} {\xi}_b^{-e_b}
\prod_{1 < c < q_3/2 \atop (c, q_3) =1} {\xi}_c^{-f_c}
\end{equation}
Since $q_1, q_2, q_3$ are mutually co-prime, 
$\Q(\zeta_{q_1}) \cap \Q(\zeta_{q_2q_3})=\Q$. So we see 
that both sides of \eqref{re} are rational numbers 
and hence equal to $\pm 1$.
Now squaring both sides, we get 
\begin{equation}
\prod_{1 < a < q_1/2 \atop (a, q_1) =1} {\xi}_a^{2d_a} 
~=~ \prod_{1 < b < q_2/2 \atop (b, q_2) =1} {\xi}_b^{-2e_b} 
\prod_{1 < c < q_3/2 \atop (c, q_3) =1} {\xi}_c^{-2f_c}~=~1.
\end{equation}
This forces that $d_a = 0$ for all
$a$ since $\xi_a$'s are multiplicatively independent. 
Again going back to
\eqref{re} and following the same argument, we get 
$e_b=0, f_c=0$ for all $b,c$.
This completes the proof.

\section{Proofs of Linear independence results}

\smallskip

\subsection{Proof of \thmref{closure}}
It is sufficient to show that given any natural number
$n$, there exist disjoint subsets
$\Omega_1, \cdots, \Omega_n \in C(q)$ 
such that $\gamma(\Omega_1, a, q), \cdots, 
\gamma(\Omega_n, a, q)$ are linearly independent over $\bQ$.
Suppose that our claim is not true. Then there exists 
an $n \in \N$ such that for
any disjoint sets $\Omega_1, \cdots, \Omega_n \in C(q)$ and
$\Omega'_1, \cdots, \Omega'_n \in C(q)$, 
we can find  $\alpha_i, \beta_j \in \bQ, ~ 1 \le i,j \le n$,
not all zero such that
$$
\alpha_1 \gamma(\Omega_1, a, q) + \cdots 
+ \alpha_n \gamma(\Omega_n, a, q) = 0 
\phantom{m}\text{ and} \phantom{m}
\beta_1 \gamma(\Omega'_1,a ,q) + \cdots 
+ \beta_n \gamma(\Omega'_n, a, q) = 0.
$$
Further assume that $\Omega_i$'s are disjoint from
$\Omega'_j$'s for all $1 \le i,j \le n$. 
Then by \eqref{br}, we have
\begin{eqnarray}\label{ga}
\gamma\sum_{i=1}^n \alpha_i \delta_{\Omega_i}
&=&
\frac{-q}{\varphi(q)} \sum_{\substack{\chi \mod q \atop \chi \ne \chi_0}}
\overline{\chi}(a) L(1,\chi) 
\sum_{i=1}^n \alpha_i \prod_{ p \in \Omega_i} ( 1 - \frac{\chi(p)}{p}) ~-~ 
\sum_{ p | q} \frac{\log p}{p-1} \sum_{i=1}^n \alpha_i \delta_{\Omega_i}  \\
\nonumber
&& \phantom{mm}
~-~~~~
\sum_{i=1}^n \alpha_i \delta_{\Omega_i} \sum_{p \in \Omega_i}
 \frac{\log p}{p-1} \\ 
\nonumber 
\end{eqnarray}
\begin{eqnarray}
\text{ and } \phantom{m}
\gamma \sum_{j=1}^n \beta_j \delta_{\Omega'_j}
&=&
\frac{-q}{\varphi(q)} \sum_{\substack{\chi \mod q \atop \chi \ne \chi_0}}
\overline{\chi}(a) L(1,\chi) 
\sum_{j=1}^n \beta_j \prod_{ p \in \Omega'_j} ( 1 - \frac{\chi(p)}{p}) ~-~ 
\sum_{ p | q} \frac{\log p}{p-1} 
\sum_{j=1}^n \beta_j \delta_{\Omega'_j} \nonumber \\
&& \phantom{mm}
~-~
\sum_{j=1}^n \beta_j \delta_{\Omega'_j} \sum_{p \in \Omega'_j}
 \frac{\log p}{p-1}. \nonumber
\end{eqnarray}
Applying \thmref{cl}, we see that $A:= 
\sum_{i=1}^n\alpha_i \delta_{\Omega_i} \ne 0$
and $B:= \sum_{j=1}^n\beta_j \delta_{\Omega'_j} \ne 0$. 
Hence from \eqref{ga}, we have
\begin{eqnarray*}
&& \frac{q}{\varphi(q)} 
\sum_{\substack{\chi \mod q \atop \chi \ne \chi_0}}
\overline{\chi}(a) L(1,\chi) 
 \left( \sum_{i=1}^n \frac{\alpha_i}{A} 
\prod_{ p \in \Omega_i} ( 1 - \frac{\chi(p)}{p}) 
~-~ 
\sum_{j=1}^n \frac{\beta_j}{B} 
\prod_{ p \in \Omega'_j} ( 1 - \frac{\chi(p)}{p}) \right) \\
&& \phantom{mmmmm} + \phantom{m}
\sum_{i=1}^n \frac{\alpha_i \delta_{\Omega_i}}{A} \sum_{p \in \Omega_i}
 \frac{\log p}{p-1}
\phantom{m} - \phantom{m}
\sum_{j=1}^n \frac{\beta_j \delta_{\Omega'_j}}{B} \sum_{p \in \Omega'_j}
 \frac{\log p}{p-1}
 ~=~ 0,
\end{eqnarray*}
a contradiction to \thmref{cl}. This completes 
the proof of \thmref{closure}.

\subsection{Proofs of \thmref{pre} and \thmref{linear-ind}}

Now we will give a proof of \thmref{pre}. 
For any finite subset $\Omega \subset P$ and
$P_{\Omega}$ as in \eqref{imp}, define
$$
S_{\Omega} := \{ u \in \N  ~|~  (u, P_{\Omega})=1 \}
$$ 
and for any natural number $u \in S_{\Omega}$, let us set
$$
\Gamma_{\Omega, u}:= \{ \gamma(\Omega, v, u) 
~|~1 \le v \le u, ~ (v, u)=1 \}.
$$
Note that the cardinality of $\Gamma_{\Omega, u}$
is $\varphi(u)$.
We claim that for any two pairwise co-prime natural 
numbers $q, r \in S_{\Omega}$, either
the set of numbers $\Gamma_{\Omega, q}$ is
linearly independent over $\Q$ or
the set of numbers $\Gamma_{\Omega, r}$
is linearly independent over $\Q$.

Suppose that our claim is not true. 
Then there exists $\alpha_a, \beta_b \in \Q$, not all zero, 
for $1 \le  a < q$ and $1 \le b < r$ 
with $(a, q)= 1= (b, r) $ such that
\begin{equation}\label{ak}
\sum_{1\le a < q \atop (a, q)=1} \alpha_a \gamma( \Omega, a, q) = 0
\phantom{m} \text{and} \phantom{m}
\sum_{1\le b < r \atop (b, r)=1} \beta_b \gamma( \Omega, b, r) = 0.
\end{equation}
Define two arithmetic functions as follows;
\begin{eqnarray}\label{panch}
f(n) 
&:=& \begin{cases}
\alpha_a  & \mbox{ if   }~ n \equiv a \!\!\!\! \pmod{q}, 
~(a, q)=1, \\
-\sum_{1 \le a \le q \atop (a, q)=1} \alpha_a 
 & \mbox{ if   }~ n \equiv 0 \!\!\!\!  \pmod{q}, \\
0  & \mbox{ otherwise},
\end{cases}\\
\phantom{m} \text{and} \phantom{m}
g(n) 
&:=& \begin{cases}
\beta_b  & \mbox{ if   }~ n \equiv b \!\!\!\! \pmod{r}, ~(b, r)=1, \\
-\sum_{1 \le b \le r \atop (b, r)=1} \beta_b 
 & \mbox{ if   }~ n \equiv 0 \!\!\!\!  \pmod{r}, \\
0  & \mbox{ otherwise}. \nonumber
\end{cases} 
\end{eqnarray}
Then $f$ and $g$ are periodic functions with periods 
$q$ and $r$ respectively. Further
$$
\sum_{1 \le a \le q} f(a) = 0
\phantom{m} \text{ and } \phantom{m}
\sum_{1 \le b \le r} g(b) = 0.
$$
Hence by \thmref{gs}, we have 
\begin{eqnarray}\label{dui}
\sum_{n \ge 1 \atop (n, P_{\Omega})=1} \frac{f(n)}{n} 
&=&  \sum_{a=1}^q f(a) \gamma( \Omega, a, q)\\
\phantom{m} \text{ and } \phantom{m}
\sum_{m \ge 1 \atop (m, P_{\Omega})=1} \frac{g(m)}{m} 
&=&  
\sum_{b=1}^r g(b) \gamma( \Omega, b, r). \nonumber 
\end{eqnarray}
Now equations \eqref{ak}, \eqref{panch}, \eqref{dui}
and the fact 
\begin{equation}\label{extra}
\gamma(\Omega, q, q) = 
\frac{1}{q}(\gamma(\Omega) - \delta_{\Omega}\log q),
\end{equation}
imply that
\begin{eqnarray*}
\sum_{n \ge 1 \atop (n, P_{\Omega})=1} \frac{f(n)}{n} 
&=&  \frac{f(q)}{q} \left( \gamma(\Omega)
 - \delta_{\Omega} \log q \right) \\
\phantom{m} \text{ and } \phantom{m}
\sum_{m \ge 1 \atop (m, P_{\Omega})=1} \frac{g(m)}{m} 
&=&  \frac{g(r)}{r} \left( \gamma(\Omega) 
- \delta_{\Omega} \log r \right).
\end{eqnarray*}
Note that $f(q)$ and $g(r)$ can not be zero.
Indeed, if for example $f(q) = 0$, then 
\begin{eqnarray}\label{char}
\sum_{n \ge 1 \atop (n, P_{\Omega})=1} \frac{f(n)}{n} 
&=& 0. 
\end{eqnarray}
However, viewing $f \chi_0$
as a periodic function modulo $qP_{\Omega}$, where 
$\chi_0$ is the trivial character modulo $P_{\Omega}$,
we have 
$$
(f\chi_0)(n) = 0 
\phantom{m} \text{ for all } \phantom{m} 
1 < (n, qP_{\Omega}) < qP_{\Omega}.
$$ 
Hence by \thmref{non-vanish}, we have
$$
\sum_{n \ge 1 \atop (n, P_{\Omega})=1} 
\frac{f(n)}{n} \ne 0,
$$
a contradiction to \eqref{char}. 
As $f(q) \ne 0$ and $g(r) \ne 0$, we have
\begin{equation}\label{tin}
\frac{q}{f(q)}\sum_{n \ge 1 \atop (n, P_{\Omega})=1} \frac{f(n)}{n} 
~+~ \delta_{\Omega} \log q
~-~ \frac{r}{g(r)} \sum_{m \ge 1 \atop (m, P_{\Omega})=1} \frac{g(m)}{m} 
~-~ \delta_{\Omega} \log r 
~=~ 0.
\end{equation}
By \thmref{gs} and equations \eqref{br}, \eqref{panch}, \eqref{extra}, 
we have
\begin{eqnarray*}
&&
\frac{q}{f(q)}\sum_{n \ge 1 \atop (n, P_{\Omega})=1} \frac{f(n)}{n} 
~+~  \delta_{\Omega} \log q \\
&=&
\frac{q}{f(q)}\sum_{1 \le a < q } 
f(a) \gamma(\Omega, a, q)
~+~ \gamma(\Omega) - \delta_{\Omega}\log q
~+~  \delta_{\Omega} \log q  \nonumber \\
\end{eqnarray*}
\begin{eqnarray*}
&=&
\frac{q}{f(q)\varphi(q)} 
\sum_{\substack{\chi \mod q \atop \chi \ne \chi_0}}
L(1,\chi)
\prod_{p\in\Omega} (1-\frac{\chi(p)}{p})
\sum_{1 \le a \le q \atop (a, qP_{\Omega})=1} f(a)\overline{\chi}(a) 
 ~~-~~ \delta_{\Omega}\sum_{p \mid q} \frac{\log p}{p-1}. 
 \end{eqnarray*}
Similarly, we have
\begin{eqnarray*}
&&\frac{r}{g(r)}\sum_{m \ge 1 \atop (m, P_{\Omega})=1} \frac{g(m)}{m} 
~+~  \delta_{\Omega} \log r \\
&=&
\frac{r}{g(r)\varphi(r)} 
\sum_{\substack{\psi \mod r \atop \psi \ne \psi_0}}
L(1,\psi)
\prod_{p\in\Omega} (1-\frac{\psi(p)}{p})
\sum_{1 \le b \le r \atop (b, rP_{\Omega})=1} g(b)\overline{\psi}(b) 
 ~~-~~ \delta_{\Omega}\sum_{p \mid r} \frac{\log p}{p-1}.
\end{eqnarray*}
Replacing these two expressions in \eqref{tin}, we see
that the left hand side of
the above expression is a non-trivial algebraic
linear combinations of $L(1, \chi)$ as $\chi$ varies over 
non-principal characters
modulo $q$, $L(1, \psi)$ as $\psi$ varies over non-principal
character modulo $r$, logarithms of prime divisors of $q$ 
and  logarithms of prime divisors of $r$. 
Then by \thmref{ESS2}, this can not be equal to zero,
a contradiction. 

Thus there exists a natural number
$r_0 \in S_{\Omega}$ such that for any $q \in \N$
with $(q, r_0P_{\Omega})=1$, the family
of numbers $\Gamma_{\Omega, q}$
are linearly independent
over $\Q$. 
Using this, we will calculate the dimension
of the space
$$
V_{\Q,N} := \Q \left <  \gamma(\Omega, m, n) ~|~   
1\le m\le n \le N \in \N, ~(m, n)=1=(n, P_{\Omega}) \right >,
$$
where $N$ is sufficiently large. To get a non-trivial lower bound on the dimension 
of $V_{\Q,N}$, we will try to find a pair of 
prime numbers $p, \ell$ in terms of $N$.

Let $t$ be the number of primes in $\Omega$. 
Now using Bertrand's Postulate, we get that there are
at least $t + 2$ primes between $\frac{N}{2^{t + 2}}$ and $N$,
where $N > 2^{t+2}$. Hence there exist
two primes $p, \ell \ge \frac{N}{2^{t+2}}$ with 
$(p\ell, P_\Omega)=1$. Thus
$$
\dim V_{\Q,N} \ge \min\{\varphi(p), \varphi(\ell)\} 
= 
\min\{p-1, \ell-1\} 
\ge \frac{N}{2^{t+2}} - 1 
\gg_\Omega N.
$$
This completes the proof of \thmref{pre}.

We now indicate the required modifications 
to derive \thmref{linear-ind}. 
Let $K$ be a number field with discriminant $d>1$
and $K \cap \Q(\zeta_{P_{\Omega}}) = \Q$.

We claim that there exists a natural number
$r_0 \in S'_{\Omega}$ such that for any $q \in \N$
with $(q, r_0 d P_{\Omega})=1$, a set of numbers
$\Gamma_{\Omega, q}$ (defined below) 
consisting of suitable $\gamma(\Omega, a, q)$'s
with $|\Gamma_{\Omega, q}| = \varphi(q)$
is linearly independent over $K$. 

In order to prove the claim, we replace 
the set $S_\Omega$ in \thmref{pre} by
$$
S'_{\Omega} := \{ u \in \N  ~|~  (u, P_{\Omega})=1, 
~ \Phi_{uP_{\Omega}}(X)
 \text{ is irreducible over } K   \}.
$$
Since $K  \cap \Q(\zeta_{uP_{\Omega}}) = \Q$ 
if and only if
$\Phi_{uP_{\Omega}}(X)$ is irreducible over $K$,
for any natural number $u$ with $(u, dP_{\Omega})=1$, 
one has $u \in S'_{\Omega}$.
Consider the set 
$$
\Gamma_{\Omega, u}:= \{ \gamma(\Omega, v, u) 
~|~1 \le v \le u, ~ (v, u)=1 \}
\text{      for   } u \in S'_{\Omega}.
$$
In order to complete the proof of the claim,
we now define $f,q$ as before.
The proof of the claim now follows
mutatis mutandis as in
\thmref{pre} except when we need to show that
neither $f(q) =0$ nor $g(r)=0$. Once again we
will use the theorem of Baker, Birch and Wirsing, 
but over number fields $K$ with $[K: \Q] > 1$. 
This forces us to have 
the additional condition that $\Phi_{uP_{\Omega}}(X)$
is irreducible over $K$. This is why we
have replaced our set $S_{\Omega}$ in 
\thmref{pre} by $S'_{\Omega}$ in \thmref{linear-ind}.

We now complete the proof of \thmref{linear-ind}. 
For a lower bound, we let
$s$ be the number of distinct prime divisors of $d$ 
and $t$ be the number of primes
in $\Omega$. Then again by Bertrand's Postulate, we get 
that there are at least  $s + t + 2$ many primes between 
$\frac{N}{2^{s+ t + 2}}$ and $N$
 with $N > 2^{s + t + 2}$.
Thus we can get two distinct primes 
$p, \ell \ge \frac{N}{2^{s + t + 2}}$
such that they are co-prime to $dP_\Omega$. Then
$$
\dim V_{K,N} \ge \min\{\varphi(p),\varphi(\ell )\} =
\min\{p-1, \ell-1\} \ge \frac{N}{2^{ s + t +2}} 
- 1 \gg_{\Omega, K} N.
$$

\subsection{Proof of \thmref{E2}}

Suppose that  $\gamma(\Omega_2, a, q), ~\gamma(\Omega_3, a, q) 
\in [\gamma(\Omega_1, a, q)]$, where 
$\Omega_1,\Omega_2$ and $\Omega_3$ are
distinct elements in $C(q)$. Then there exist non-zero algebraic numbers
$\beta, \lambda$ such that
\begin{equation}\label{rel1}
\gamma(\Omega_1, a, q) ~=~ \beta \gamma(\Omega_2, a, q) 
\phantom{m} 
\text{and}
\phantom{m}
\gamma(\Omega_1, a, q) ~=~  \lambda \gamma(\Omega_3, a, q).
\end{equation}
For a Dirichlet character $\chi$ modulo $q$ and 
a finite set $\Omega$ consisting of primes
co-prime to $q$, we define 
\begin{eqnarray*}
a_\Omega 
&:=&  
\frac{\delta_\Omega}{q}  \ne 0, \phantom{m}
\gamma_1 
~:=~
\gamma + \sum_{p \mid q} \frac{\log p}{p-1} \\
\text{ and} \phantom{m}
\alpha_{\chi,\Omega} 
&:=&
\frac{\overline{\chi}(a)}{\varphi(q)} 
\prod_{p\in\Omega} (1-\frac{\chi(p)}{p}).
\end{eqnarray*}
Using \eqref{br} and \eqref{rel1}, we get
\begin{eqnarray}\label{eq61}
&& 
\gamma_1(a_{\Omega_1}-\beta a_{\Omega_2})
+
a_{\Omega_1}\sum_{\substack{p\in\Omega_1}}\frac{\log p}{p-1}
-
\beta a_{\Omega_2}\sum_{\substack{p\in\Omega_2}} \frac{\log p}{p-1} 
+
\sum_{\substack{\chi \mod q \atop \chi \ne \chi_0}}L(1,\chi)
(\alpha_{\chi,\Omega_1}-\beta\alpha_{\chi,\Omega_2}) 
~=~ 0.
\end{eqnarray}
Similarly, we have
\begin{eqnarray}\label{eq71}
&&
\gamma_1(a_{\Omega_1}-\lambda a_{\Omega_3})
+ 
a_{\Omega_1}\sum_{\substack{p\in\Omega_1}}\frac{\log p}{p-1}
-
\lambda a_{\Omega_3}\sum_{\substack{p\in\Omega_3}} 
\frac{\log p}{p-1} 
+ 
\sum_{\substack{\chi \mod q \atop \chi \ne \chi_0}}L(1,\chi)
(\alpha_{\chi,\Omega_1} 
- 
\lambda\alpha_{\chi,\Omega_3})
~=~ 0. 
\end{eqnarray}
Since $\Omega_1,\Omega_2,\Omega_3$ are distinct set of primes,
applying \thmref{cl} to equations \eqref{eq61} and \eqref{eq71}, we get 
$$
a_{\Omega_1} - \beta a_{\Omega_2} \ne 0, 
\phantom{m}
a_{\Omega_1} - \lambda a_{\Omega_3} \ne 0.
$$
By the same reasoning, we see that
$$
a_{\Omega_3}-\frac{\beta}{\lambda} a_{\Omega_2}\ne 0.
$$
Again from \eqref{eq61} and \eqref{eq71}, it follows that 
\begin{eqnarray*}
\frac{a_{\Omega_1}(\beta a_{\Omega_2}-\lambda a_{\Omega_3})}
{(a_{\Omega_1}-\beta a_{\Omega_2})(a_{\Omega_1}-\lambda a_{\Omega_3})}
\sum_{p\in \Omega_1}\frac{\log p}{p-1} 
~-~ \frac{\beta a_{\Omega_2}}{(a_{\Omega_1}-\beta a_{\Omega_2})} 
\sum_{p\in \Omega_2}\frac{\log p}{p-1}\\ 
~+~\frac{\lambda a_{\Omega_3}}{(a_{\Omega_1}-\lambda a_{\Omega_3})}
\sum_{p\in \Omega_3}\frac{\log p}{p-1}
~+~\sum_{\substack{\chi \mod q 
\atop \chi \ne \chi_0}}L(1,\chi)A(\chi)
~=~ 0,
\end{eqnarray*}
where
$$
A(\chi)=\frac{(\alpha_{\chi,\Omega_1}-\beta \alpha_{\chi,\Omega_2})}
{(a_{\Omega_1}-\beta a_{\Omega_2})}~-~
\frac{(\alpha_{\chi,\Omega_1}-
\lambda\alpha_{\chi,\Omega_3})}{(a_{\Omega_1}-\lambda a_{\Omega_3})}.
$$
Since $\Omega_1,\Omega_2,\Omega_3$ are distinct set of primes, 
without loss of generality, one can assume that there exists a prime
$p_1 \in \Omega_1$ such that either $p_1 \notin \Omega_2 
\cup \Omega_3$ or $p_1 \in \Omega_2$ but not in 
$\Omega_3$. The coefficient of $\log p_1$ in the first case is 
$$
\frac{a_{\Omega_1}(\beta a_{\Omega_2}-\lambda a_{\Omega_3})}
{(a_{\Omega_1}-\beta a_{\Omega_2})(a_{\Omega_1}
-\lambda a_{\Omega_3})(p_1-1)} \ne 0
$$
and in the second case is 
$$
\frac{\lambda a_{\Omega_3}}{(\lambda a_{\Omega_3}
- a_{\Omega_1})(p_1-1)} \ne 0.
$$
Hence in both cases we arrive at a contradiction by \thmref{cl}.

\subsection{Proof of \thmref{E3}}

Suppose that $\gamma(\Omega,a,q_2), ~\gamma(\Omega,a,q_3) 
\in [\gamma(\Omega, a, q_1)]$, where $q_1, q_2, q_3$
are distinct elements in $C(a, \Omega)$.
Then there exist non-zero algebraic numbers
$\beta, \lambda$ such that
\begin{equation}\label{rel}
\gamma(\Omega, a, q_1) ~=~ \beta \gamma(\Omega, a, q_2),
\phantom{m}
\text{and}
\phantom{m}
\gamma(\Omega, a, q_1) ~=~  \lambda \gamma(\Omega, a, q_3).
\end{equation}
Write 
\begin{eqnarray*}
a_{q_i}  &:=&  
\frac{\delta_\Omega}{q_i} \ne 0, \phantom{m}
\gamma_1 
~:=~
\gamma + \sum_{p \mid \Omega} \frac{\log p}{p-1} \\
\text{ and} \phantom{m}
\alpha_{\chi, q_i} 
&:=&
\frac{\overline{\chi}(a)}{\varphi(q_i)} 
\prod_{p\in\Omega} (1-\frac{\chi(p)}{p}).
\end{eqnarray*}
Using \eqref{br} and \eqref{rel}, we get
\begin{eqnarray}\label{eq6}
\gamma_1(a_{q_1} - \beta a_{q_2})
& + &
a_{q_1}\sum_{\substack{p \mid q_1}}\frac{\log p}{p-1}
-
\beta a_{q_2}\sum_{\substack{p \mid q_2}} \frac{\log p}{p-1}  
+
\sum_{\substack{\chi \mod q_1 \atop \chi \ne \chi_0}}
\alpha_{\chi, q_1} L(1,\chi) \\
&&  
\phantom{mm}
-  \beta \sum_{\substack{\chi \mod q_2 \atop \chi \ne \chi_0}}
\alpha_{\chi, q_2}L(1, \chi)
~=~ 0. \nonumber
\end{eqnarray}
Similarly, we have
\begin{eqnarray}\label{eq7}
\gamma_1(a_{q_1} - \lambda a_{q_3})
&+&
a_{q_1}\sum_{\substack{p \mid q_1}}\frac{\log p}{p-1}
-
\lambda a_{q_3}\sum_{\substack{p \mid q_3}} \frac{\log p}{p-1}  
+
\sum_{\substack{\chi \mod q_1 \atop \chi \ne \chi_0}}
\alpha_{\chi, q_1} L(1,\chi) \\
&&  
\phantom{mm}
-  \lambda \sum_{\substack{\chi \mod q_3 \atop \chi \ne \chi_0}}
\alpha_{\chi, q_3}L(1, \chi)
~=~ 0.  \nonumber 
\end{eqnarray}
Since $q_1, q_2$ and $q_3$ are mutually co-prime natural numbers,
applying \thmref{ESS2} to equations \eqref{eq6} and \eqref{eq7}, we get 
$$
a_{q_1} - \beta a_{q_2} \ne 0, 
\phantom{m}
a_{q_1} - \lambda a_{q_3} \ne 0.
$$
Similar reasoning show that
$$
\beta a_{q_2}-\lambda a_{q_3} \ne 0.
$$
Hence we have
\begin{eqnarray*}
C a_{q_1}
\sum_{p\mid q_1}\frac{\log p}{p-1} 
~-~
 \frac{\beta a_{ q_2}}{(a_{q_1}-\beta a_{q_2})} 
\sum_{p\mid q_2}\frac{\log p}{p-1}
~+~
\frac{\lambda a_{q_3}}{(a_{q_1} - \lambda a_{q_3})}
\sum_{p\mid q_3}\frac{\log p}{p-1} 
~+~
C \sum_{\substack{\chi \mod q_1 \atop \chi \ne \chi_0}}
\alpha_{\chi, q_1} L(1,\chi)\\
~+~ 
\frac{\lambda}{(a_{q_1} - \lambda a_{q_3})} 
\sum_{\substack{\chi \mod q_3 \atop \chi \ne \chi_0}}
\alpha_{\chi, q_3} L(1, \chi) 
~-~ \frac{\beta }{(a_{q_1} - \beta a_{q_2})}
\sum_{\substack{\chi \mod q_2 \atop \chi \ne \chi_0}}
\alpha_{\chi, q_2}L(1, \chi) 
~=~ 0,
\end{eqnarray*}
where
$$
C:= \frac{\beta a_{q_2}-\lambda a_{q_3}}
{(a_{q_1}-\beta a_{q_2})(a_{q_1}-\lambda a_{q_3})} \ne 0,
$$
a contradiction to \thmref{ESS2}. 
This completes the proof of the theorem.

\smallskip

\section{Consequences of weak Schanuel conjecture}

\smallskip

In this section, we state some conditional results on algebraic
independence of Euler-Briggs constants assuming the
Weak Schanuel conjecture (see page 111 of \cite{MR}, see also \cite{GMR}). 

\begin{conj}{\rm(Weak Schanuel)}~\label{S}
Let $\alpha_1, \cdots, \alpha_n$ be non-zero algebraic 
numbers such that the numbers $\log \alpha_1, \cdots, \log \alpha_n$
are $\Q$-linearly independent. Then 
$\log \alpha_1, \cdots, \log \alpha_n$
are algebraically independent.
\end{conj}

Before we proceed further, let us fix few more notation.
For $a, q \in \N$ with $1 \le a \le q$ and $(a, q)=1$,
we define 
$$
\gamma^*(\Omega, a, q) := \frac{q\gamma(\Omega, a, q)}{\delta_{\Omega}}.
$$
We say $\gamma^*(\Omega_1, a, q) \sim \gamma^*(\Omega_2, a, q)$ if
there exists a non-zero algebraic number $\alpha$ such that
$$
\gamma^*(\Omega_1, a ,q)  = \alpha \gamma^*(\Omega_2, a, q).
$$
Note that $\gamma(\Omega_1, a, q) \sim \gamma(\Omega_2, a ,q)$ if and only if
$\gamma^*(\Omega_1, a, q) \sim \gamma^*(\Omega_2, a, q)$. Hence we will
study $\gamma^*(\Omega, a, q)$ in place of $\gamma(\Omega, a, q)$
whenever convenient.

We call a finite sequence of sets
$\{ \Omega_1, \cdots , \Omega_n \}$ 
an irreducible  sequence if 
$$
\cup_{i=1}^n \Omega_i  \ne \cup_{j \in J} \Omega_j
$$
for any proper subset $J \subset \{ 1, \cdots, n\}$.
We call an infinite sequence of distinct sets 
$\{ \Omega_n ~|~  n \in \N\}$ an irreducible sequence
if every finite subsequence is irreducible. 
It is easy to see if 
$$
p_1 < p_2< \cdots 
$$ 
is a sequence of distinct prime numbers
and $\Omega_i = \{ p_i  \}$, then $\{ \Omega_i \}$
is an irreducible sequence.  On the other hand,
the sequence 
$$
 \{ p_1 \}, \{p_2 \},  \{p_1, p_2 \}, \{ p_3 \},
\{p_1, p_2, p_3 \}, \cdots 
$$ 
where $p_i$'s are distinct
prime numbers is not an irreducible sequence
though it contains an irreducible
subsequence.
Here we have the following theorem.

\begin{thm}\label{E5}
Suppose that the Weak Schanuel conjecture is true. 
Further, suppose that
$T := \{ \Omega_n  \}_{n \in \N} $ is an
infinite sequence of non-empty finite subsets 
of prime numbers co-prime to $q$.
Consider the set  
$$
S_1:= \{ \gamma^*(\Omega_n, a, q) - \gamma 
- \sum_{\chi \ne \chi_0} \alpha^*_{\chi, \Omega_n, q} L(1, \chi)
~ | ~ \Omega_n \in T \},
$$ 
where $\chi$ is a Dirichlet character modulo $q$ and
$$
\alpha^*_{\chi,\Omega, q}:= \overline{\chi}(a)
\prod_{p\in\Omega} (1-\frac{\chi(p)}{p})(1-\frac{1}{p})^{-1}
\prod_{p | q} (1-\frac{1}{p})^{-1}.
$$
Then the elements of $S_1$ are algebraically independent 
if the infinite sequence $T$ is irreducible. 
\end{thm}

\smallskip
\noindent
{\bf Proof of \thmref{E5}.}
Let $T = \{ \Omega_n \}_{n \in \N}$ be
an irreducible sequence where no $\Omega_n$
contains any prime divisors of $q$. 
By \eqref{br}, we know that
$$
A_n := \gamma^*(\Omega_n, a, q) - \gamma 
- \sum_{\chi \ne \chi_0} \alpha^*_{\chi, \Omega_n, q} ~L(1, \chi)
= \sum_{p \in \Omega_n}\frac{\log p}{p-1} 
+ \sum_{p \mid q}\frac{\log p}{p-1}.
$$
Hence by Weak Schanuel's conjecture \ref{S}, it is sufficient
to show that the elements $A_n$'s
for $\Omega_n \in T$ are linearly independent over $\Q$.
If not, then there exists a finite subsequence
$T' = \{ \Omega_{n_1}, \cdots, \Omega_{n_k} \}$ 
of $T$ and integers $m_1, \cdots, m_k$, not all zero, 
such that
\begin{equation}\label{eq4}
m_1 A_{n_1}  ~+~ \cdots  ~+~  m_k A_{n_k}  = 0
\end{equation}
Write $ \Omega := \cup_{i=1}^k \Omega_{n_i}$.
Then applying \eqref{br} in \eqref{eq4}, we get 
\begin{equation}\label{new}
\sum_{ p \in \Omega } t_p \log p
~+~   \sum_{ \ell | q} r_{\ell} \log \ell 
~=~ 0,
\end{equation}
where $t_p, r_{\ell} \in \Q$ and $p \in \Omega$ with $(p,q)=1$.
Since $T'$ is an irreducible sequence and not all 
$m_i$'s are zero, it follows that not all $t_p$'s are zero,
a contradiction to \eqref{new}.
This completes the proof of \thmref{E5}. 

\bigskip

Before we state our next theorem, let us introduce a notation and a definition. 
For $I \subseteq \N$, let $P(I)$ be the set 
of all prime divisors of the elements of $I$. A finite subset $I$ of $\N$
is called irreducible if and only if
\begin{equation}\label{irr}
P(I)\neq\cup_{J\subsetneq I}P(J).
\end{equation}
An infinite subset $T\subseteq\N$ is called irreducible if all
finite subset of $T$ are irreducible. In this context, we 
have the following theorem.

\begin{thm}\label{E6}
Suppose that the Weak Schanuel conjecture is true. Let $\Omega$ be a
finite set of primes. Further, suppose that $T=\{q_i\}_{i\in\N}$
be an infinite irreducible sequence of natural numbers co-prime to the primes
in $\Omega$. Let $a\in\N$ be such that $(a,q_i)=1$ for all $i\in\N$. 
Consider the set
$$
S_2:= \{ \gamma^*(\Omega, a, q) - \gamma 
- \sum_{\chi \ne \chi_0} \alpha^*_{\chi, \Omega, q} L(1, \chi)
~ | ~ q\in T \},
$$ 
where $\alpha^*_{\chi, \Omega, q}$ is as in \thmref{E5}.
Then the elements of $S_2$ are algebraically independent.
\end{thm}

\smallskip
\noindent
{\bf Proof of \thmref{E6}.}
We have
$$
\gamma^*(\Omega, a, q) - \gamma 
- \sum_{\chi \ne \chi_0} \alpha^*_{\chi,\Omega, q} L(1, \chi)
= \sum_{p \in \Omega}\frac{\log p}{p-1} + \sum_{p \mid q}\frac{\log p}{p-1}.
$$
Hence by Weak Schanuel's conjecture \ref{S}, it is sufficient 
to show that the elements of $S_2$ are 
linearly independent over $\Q$.
Suppose not, then there exists a finite subset $\{q_1,\cdots,q_n\}$
of $T$ and integers $m_1,\cdots,m_n$, not all $0$, such that
$$
\sum_{i=1}^{n}m_i\{\gamma^*(\Omega, a, q_i) - \gamma 
- \sum_{\chi \ne \chi_0} \alpha^*_{\chi, \Omega, q_i} L(1, \chi)\}=0.
$$
So
$$
\sum_{p \in \Omega}\frac{\log p}{p-1}\sum_{i=1}^{n}m_i
~~+~~ \sum_{i=1}^{n}m_i\sum_{p \mid q_i}\frac{\log p}{p-1}=0.
$$
Without loss of generality let $m_1$ be non-zero. 
Since $T$ is irreducible, by definition 
all finite subsets of $T$ are  irreducible.
Then by \eqref{irr} we know that
$P(q_1, q_2, \cdots,q_n) \ne P(q_2, q_3\cdots,q_n)$
 and hence there exists a prime $p$ such that $p\mid q_1$ 
but $p\nmid q_j$ for all $j\neq 1$. This implies that 
the coefficient of $\log p$ is
$m_1/(p-1) \neq 0$, a contradiction by Baker's theorem.

\medskip

\bigskip
\noindent
{\bf Acknowledgments.}
Part of the work of the first and third author was
supported by a DAE number theory grant. 
The work of the first author was also 
supported by her SERB grant. Some part of the work was 
done when the first author was visiting 
the Mathematics department of
University of Bordeaux as an ALGANT 
scholar. She would like to thank the members
of the institute for their hospitality and the ALGANT
program for the financial support. We would also like
to thank the referee for his/her suggestions
which improved the content as well as presentation
of the paper.

\end{document}